\documentclass[reqno]{amsart}
\usepackage{amssymb,amsmath}
\usepackage{verbatim}
\usepackage{mathrsfs}
\usepackage[all]{xy}
\usepackage[pdftex]{graphicx}
\usepackage[pdftex]{hyperref}

\title{The Non-Archimedean Monge-Amp\`ere Equation}
\author{S{\'e}bastien Boucksom
  \and
  Charles Favre
  \and
  Mattias Jonsson}
\address{CNRS-CMLS\\
  \'Ecole Polytechnique\\
  F-91128 Palaiseau Cedex\\
  France}
\email{boucksom@polytechnique.edu}
\address{CNRS-CMLS\\
  \'Ecole Polytechnique\\
  F-91128 Palaiseau Cedex\\
  France}
\email{charles.favre@polytechnique.edu}
\address{Dept of Mathematics\\
  University of Michigan\\
  Ann Arbor, MI 48109-1043\\
  USA}
\email{mattiasj@umich.edu}

\date\today

\numberwithin{equation}{section}       

\theoremstyle{plain}
\newtheorem{Thm}{Theorem}[section]

\newtheorem{Prop-def}[Thm]{Proposition-Definition}

\theoremstyle{definition}

\theoremstyle{remark}

\newcounter{alphasect}
\def\alphainsection{0}

\let\oldsection=\section
\def\section{%
  \ifnum\alphainsection=1%
\e    \addtocounter{alphasect}{1}
  \fi%
\oldsection}%

\renewcommand\thesection{%
 \ifnum\alphainsection=1%
   \Alph{alphasect}%
 \else
   \arabic{section}%
 \fi%
}%

\newcommand{\C}{{\mathbf{C}}}

\newcommand{\G}{{\mathbf{G}}}

\newcommand{\Q}{{\mathbf{Q}}}
\newcommand{\R}{{\mathbf{R}}}
\newcommand{\Z}{{\mathbf{Z}}}

\newcommand{\cE}{{\mathcal{E}}}

\newcommand{\cH}{{\mathcal{H}}}

\newcommand{\cL}{{\mathcal{L}}}

\newcommand{\cO}{{\mathcal{O}}}

\newcommand{\cX}{{\mathcal{X}}}

\newcommand{\fa}{{\mathfrak{a}}}

\renewcommand{\a}{\alpha}

\newcommand{\e}{\varepsilon}
\newcommand{\f}{\varphi}

\newcommand{\self}{\circlearrowleft}

\newcommand{\ie}{i.e.\ }

\newcommand{\Lan}{{L^\mathrm{an}}}
\newcommand{\Tan}{{T^\mathrm{an}}}

\newcommand{\Xan}{{X^\mathrm{an}}}

\newcommand{\Xtrop}{{X^\mathrm{trop}}}

\newcommand{\can}{\operatorname{can}}

\renewcommand{\div}{\operatorname{div}}

\newcommand{\MA}{\operatorname{MA}}

\newcommand{\charac}{\operatorname{char}}

\newcommand{\PSH}{\operatorname{PSH}}

\newcommand{\Spec}{\operatorname{Spec}}

\newcommand{\trop}{\operatorname{trop}}

\newcommand{\Vol}{\operatorname{Vol}}

\newcommand{\cro}[1]{[\![#1]\!]}
\newcommand{\lau}[1]{(\!(#1)\!)}

\begin{document}
\begin{abstract}
  We give an introduction to our work on the 
  solution to the non-Archimedean Monge-Amp\`ere equation
  and make comparisons to the complex counterpart.
  These notes are partially based on talks at the 
  2015 Simons Symposium on Tropical and Nonarchimedean Geometry.
\end{abstract}

\maketitle
%
%
%
%
%
%
\section*{Introduction}
The purpose of these notes is to discuss the Monge-Amp\`ere
equation
\begin{equation*}
  \MA(\phi)=\mu
\end{equation*}
in both the complex and non-Archimedean setting. 
Here $\mu$ is a positive
measure\footnote{All measures in this paper will be assumed to be Radon measures.} 
on the analytification of a smooth projective variety,
$\phi$ is a semipositive metric on an ample line bundle on $X$, 
and $\MA$ is the Monge-Amp\`ere operator.
All these terms will be explained below.

In the non-Archimedean case, our presentation is based on the 
papers~\cite{siminag,nama} to which we refer for details. 
In the complex case, we follow~\cite{BBGZ} rather closely.
Generally speaking, we avoid technicalities or detailed proofs.

\medskip\noindent
\textbf{Acknowledgement}.
These notes are partially based on talks given at the 
2015 Simons Symposium on Tropical and Nonarchimedean Geometry
by the first and third authors.
We want to thank the Simons foundation as well as the organizers
and participants of the symposium for the opportunity to
present our work, and for stimulating discussions. 
Boucksom was supported by the ANR projects MACK and POSITIVE\@. 
Favre was supported by the ERC-starting grant project ”Nonarcomp” no.307856.
Jonsson was supported by NSF grant DMS-1266207. 
%
%
%
%
%
%
\section{Metrics on lines bundles}\label{S101}
Let $K$ be a field equipped with a complete multiplicative norm
and let $X$ be a smooth projective variety over $K$.
To this data we can associate an analytification $\Xan$. 
When $K$ is the field of complex 
numbers with its usual norm, $\Xan$ is a compact complex manifold.
When the norm is non-Archimedean, $\Xan$ is a $K$-analytic
space in the sense of Berkovich~\cite{BerkBook}.
In either case, it is a compact Hausdorff space.

Let $L$ be a line bundle on $X$. It also admits an analytification $\Lan$.
A \emph{metric} on $\Lan$ is
a rule that to a local section 
$s:U\to\Lan$, where $U\subset\Xan$, associates a 
function $\|s\|$ on $U$,
subject to the condition $\|fs\|=|f|\cdot\|s\|$,
for any analytic function $f$ on $U$. The metric
is \emph{continuous} if $\|s\|$ is continuous
on $U$ for every $s$.

For our purposes it is convenient to use additive notation for
metrics and line bundles. Given an open cover $U_\a$ of $\Xan$ and local
trivializations of $\Lan$ on each $U_\a$, we can identify
a section $s$ of $L$ with a collection $(s_\a)_\a$ of analytic 
functions. A metric $\phi$ is then a collection of functions 
$(\phi_\a)_\a$ in such a way that $\|s\|_\phi=|s_\a|e^{-\phi_\a}$ on $U_\a$.
With this convention, if $\phi$ is a metric on $\Lan$, any other
metric is of the form $\phi+f$, where $f$ is a function on $\Xan$.
If $\phi_i$ is a metric on $L_i$, $i=1,2$, then $\phi_1+\phi_2$
is a metric on $L_1+L_2$.

Over the complex numbers, smooth metrics $\phi$ 
(\ie each $\phi_\a$ is smooth), play an important role. 
Of similar status, for $K$ non-Archimedean, are \emph{model metrics}
defined as follows.\footnote{Model metrics are not smooth in the sense
  of~\cite{CD12} but nevertheless, for our purposes, play the same
  role as smooth metrics in the complex case.}
Let $R$ be the valuation ring of $K$ and $k$ the residue field.
A \emph{model} of $X$ is a normal scheme $\cX$, flat
and projective over $\Spec R$ and with generic fiber isomorphic
to $X$. A model of $L$ is a $\Q$-line bundle $\cL$ on
$\cX$ whose restriction to $X$ is isomorphic to $L$. It defines a
continuous metric $\phi_\cL$ on $L$ in such a way that any local 
nonvanishing section of a muliple of $\cL$ has norm constantly equal to one.
Model functions, that is, model metrics on $\cO_X$, are  
dense in $C^0(\Xan)$.  
 We refer to~\cite{Ch2} or~\cite{siminag} for a more thorough discussion.

Over $\C$, a smooth metric $\phi$ on $\Lan$ is semipositive (positive) if its
curvature form $dd^c\phi$ is a semipositive (positive) $(1,1)$-form.
Here
$dd^c\phi=dd^c\phi_\a=\frac{i}{\pi}\partial\overline{\partial}\phi_\a$ 
for any $\a$. Such metrics only exist when $L$ is nef.

In the non-Archimedean setting we say that a model metric $\phi_\cL$
on $\Lan$ is semipositive if the line bundle $\cL$ is relatively nef,
that is, its degree is nonnegative 
on any proper curve contained in the special fiber $\cX_0$.
This implies that $L$ is nef. 

In both the complex and non-Archimedean case we say that a 
continuous metric $\phi$ is semipositive if there exists  
a sequence $(\phi_m)_1^\infty$ of semipositive smooth/model metrics 
such that $\lim_{m\to\infty}\sup_{\Xan}|\phi_m-\phi|=0$.  
In the non-Archimedean case, this notion was first introduced 
by Zhang~\cite{Zha95a} and Gubler~\cite{Gub98}.
In the complex case, it is more natural to say that a 
continuous metric $\phi$ is semipositive if its curvature current 
$dd^c\phi$ is a positive closed current. At least when $L$ is ample,
one can then prove 
(see~\S\ref{S106} below) that $\phi$ can be approximated by 
smooth metrics; such an approximation is furthermore crucial 
for many arguments in pluripotential theory.

In the non-Archimedean case, Chambert-Loir and Ducros have 
introduced a notion of forms and currents on Berkovich spaces. However,
it is not known whether a continuous metric whose curvature
current (in their sense) is semipositive can be approximated
by semipositive model metrics.

In both the complex and non-Archimedean case we denote by 
$\PSH^0(\Lan)$ the space of continuous semipositive metrics 
on $\Lan$. Here the superscript refers to continuity ($C^0$) whereas 
``PSH'' reflects the fact that in the complex case, semipositive metrics 
are global versions of plurisubharmonic functions.
%
%
%
%
%
%
\section{The Monge-Amp\`ere operator}\label{S102}
In the complex case, the Monge-Amp\`ere operator is a 
second order differential operator:
we set $\MA(\phi)=(dd^c\phi)^n$ for a smooth metric $\phi$.
It is a nonlinear operator if $n>1$.
When $\phi$ is semipositive, $\MA(\phi)$ is a smooth 
positive measure on $\Xan$ of mass $(L^n)$. 
It is a \emph{volume form}, that is, 
equivalent to Lebesgue measure, if $\phi$ is positive. 

Next we turn to the non-Archimedean setting. 
\emph{From now on we assume that $K$ is discretely valued}. 
Pick a uniformizer $t$ of the maximal ideal in the
valuation ring $R$ of $K$. 

Consider a model metric $\phi_\cL$, associated
to a model $(\cX,\cL)$ of $(X,L)$ over $\Spec R$.
Write the special fiber as $\cX_0=\div t=\sum_{i\in I}b_iE_i$,
where $E_i$ are the irreducible components of $\cX_0$ and
$b_i\in\Z_{>0}$. To each $E_i$
is associated a unique (divisorial) point $x_i\in\Xan$.
We then define 
\begin{equation*}
  \MA(\phi):=\sum_{i\in I}b_i(\cL|_{E_i})^n\delta_{x_i}.
\end{equation*}
If $\phi_\cL$ is semipositive, $\cL|_{E_i}$ is nef; hence 
$(\cL|_{E_i})^n\ge0$ and $\MA(\phi_\cL)$ is a positive measure.
Its total mass is
\begin{equation*}
  \int_{\Xan}1\cdot\MA(\phi_\cL)
  =\sum_{i\in I}b_i(\cL|_{E_i})^n
  =(\cL^n\cdot\cX_0)
  =(\cL^n\cdot\cX_\eta)
  =(L^n).
\end{equation*}
Here the second to last equality follows from the flatness of $\cX$ 
over $\Spec R$,
and the last equality from $\cX_\eta\simeq X$, $\cL_\eta\simeq L$.

From now on assume that $L$ is \emph{ample}, that is, we have a
\emph{polarized} pair $(X,L)$.
In both the complex and non-Archimedean case we define 
$\MA(\phi)$ for a continuous semipositive metric by
$\MA(\phi):=\lim_{m\to\infty}\MA(\phi_m)$ for any sequence $(\phi_m)_1^\infty$
converging uniformly to $\phi$. Of course, it is not obvious
that the limit exists or independent of the sequence $(\phi_m)_1^\infty$.
In the complex case this is a very special case of the Bedford-Taylor theory
developed in~\cite{BT1,BT2}. The analogous analysis in the 
non-Archimedean case is due to Chambert-Loir~\cite{Ch1}.
%
%
%
%
%
%
\section{The complex Monge-Amp\`ere equation}\label{S103}
\begin{Thm}
  Let $(X,L)$ be a polarized complex projective variety of dimension $n$
  and let $\mu$ be a positive measure on $\Xan$ of total mass $(L^n)$.
  \begin{itemize}
  \item[(i)]
    If $\mu$ is a volume form, then there exists a smooth positive
    metric $\phi$ on $\Lan$ such that $\MA(\phi)=\mu$. 
  \item[(ii)]
    If $\mu$ is absolutely continuous with respect to Lebesgue measure, with 
    density in $L^p$ for some $p>1$, then there
    exists a (H\"older) continuous metric $\phi$ on $\Lan$ 
    such that $\MA(\phi)=\mu$. 
  \item[(iii)]
    The metrics in~(i) and~(ii) are unique up to additive constants.
  \end{itemize}
\end{Thm}
The uniqueness statement in the setting of~(i) is due to Calabi. 
The much harder existence part was proved by Yau~\cite{Yau78}, 
using PDE techniques. The combined result is often called the 
Calabi-Yau Theorem.

The general setting of~(ii)--(iii) was treated by
Ko{\l}odziej~\cite{kolodziej1,kolodziej2} who 
used methods of pluripotential theory together with 
a nontrivial reduction to Yau's result.
Guedj and Zeriahi~\cite{GZ2} more generally established the existence of 
solutions of $\MA(\phi)=\mu$ for positive measures $\mu$ 
(of mass $(L^n)$) that do not put mass on pluripolar sets.
In this generality,
the metrics $\phi$ are no longer continuous but rather lie in a 
suitable energy class, modeled upon work by Cegrell~\cite{Ceg2}.
Dinew~\cite{Dinew}, improving upon an earlier result by
B{\l}ocki~\cite{Blocki}, proved the corresponding uniqueness theorem. 
All these existence and uniqueness results are furthermore valid 
(in a suitable formulation) in the 
transcendental case, when $(X,\omega)$ is a K\"ahler manifold.

The complex Monge-Amp\`ere equation is of fundamental importance
to complex geometry. For example, it implies that every compact complex 
manifold with vanishing first Chern class (such manifolds are now called
Calabi-Yau manifolds) admit a Ricci flat metric in any given K\"ahler class.
The complex Monge-Amp\`ere equation also plays a key role in recent work
on the space of K\"ahler metrics.
%
%
%
%
%
%
\section{The non-Archimedean Monge-Amp\`ere equation}\label{S104}
As before, suppose $K\simeq k\lau{t}$ is a discretely valued field 
with valuation ring $R\simeq k\cro{t}$
and residue field $k$. We further assume that $K$ has 
\emph{residue characteristic zero}, $\charac k=0$.
This implies that $R\simeq k\cro{t}$ and $K\simeq k\lau{t}$, where $k$
is the residue field of $K$.
More importantly, 
$X$ then admits SNC models, that is, regular models $\cX$ such that
the special fiber $\cX_0$ has simple normal crossings. The dual complex 
$\Delta_\cX$, encoding intersections between irreducible components of $\cX_0$,
then embeds as a compact subset of $\Xan$.
\begin{Thm}\label{T103}
  Let $(X,L)$ be a polarized complex projective variety of dimension $n$
  over $K$. Assume $X$ is defined over a smooth $k$-curve.
  Let $\mu$ be a positive measure on $\Xan$ of total mass $(L^n)$,
  supported on the dual complex of some SNC model.
  \begin{itemize}
  \item[(i)]
    There exists a continuous metric $\phi$ on $\Lan$ such that $\MA(\phi)=\mu$. 
  \item[(ii)]
    The metric in~(i) is unique up to an additive constant.
  \end{itemize}
\end{Thm}
Here the condition on $X$ means that there exists a smooth projective 
curve $C$ over $k$, a smooth projective variety $Y$ over $C$, and a point 
$p\in C$ such that $X$ is isomorphic to the base change 
$Y\times_k\Spec K$, where $K$ is the fraction field of $\widehat{\cO}_{C,p}$.
This condition is presumably redundant, but is used in the 
proof: see~\S\ref{S108}.

To our knowledge, the first to consider the Monge-Amp\`ere equation
(or Calabi-Yau problem) in a non-Archimedean setting were
Kontsevich and Tschinkel~\cite{KoTs}. They outlined a strategy in the case
when $\mu$ is a point mass. 

The case of curves ($n=1$) was
treated in detail by Thuillier in his thesis~\cite{thuillierthesis};
see also~\cite{BRbook,valtree}. In this case,
the Monge-Amp\`ere equation is linear and one can construct fundamental
solutions by exploring the topological structure of $\Xan$.

In higher dimensions, Yuan and Zhang~\cite{yuanzhang} proved the uniqueness
statement~(ii). Their proof, based on the method by B{\l}ocki, 
is valid in a more general context than stated above.
The first existence result was obtained by Liu~\cite{Liu11},
who treated the case when $X$ is a maximally degenerate abelian variety and 
$\mu$ is equivalent to Lebesgue measure on the skeleton of $X$. 
His approach amounts to solving a \emph{real} Monge-Amp\`ere equation
on the skeleton.
The existence result~(i) above was proved by the authors 
in~\cite{nama} and the companion paper~\cite{siminag}.
We will discuss our approach below.

The geometric ramifications of the non-Archimedean 
Monge-Amp\`ere equations remain to be developed.
%
%
%
%
%
%
\section{A variational approach}\label{S105}
We shall present a unified approach to solving the complex and 
non-Archimedean Monge-Amp\`ere equations in any dimension.
The method goes back to Alexandrov's work in convex 
geometry~\cite{Alexandrov}. It was adapted to the complex case
in~\cite{BBGZ} and to the non-Archimedean analogue in~\cite{nama}.

The general strategy is to construct an \emph{energy functional}
\begin{equation*}
  E:\PSH^0(\Lan)\to\R
\end{equation*}
whose derivative is the Monge-Amp\`ere operator, $E'=\MA$,
in the sense that 
\begin{equation*}
  \frac{d}{dt}E(\phi+tf)|_{t=0}=\int_{\Xan}f\MA(\phi),
\end{equation*}
for every continuous semipositive metric 
$\phi\in\PSH^0(\Lan)$ and every smooth/model function $f$
on $\Xan$.

Grant the existence of this functional for the moment.
Given a measure $\mu$ on $\Xan$,  consider the functional 
$F_\mu:\PSH^0(\Lan)\to\R$ defined by 
\begin{equation*}\
  F_\mu(\phi)=E(\phi)-\int\phi\,\mu.
\end{equation*}
Suppose we can find $\phi\in\PSH^0(\Lan)$ that maximizes $F_\mu$.
Since the derivative of $F_\mu$ is equal to $F_\mu'=\MA-\mu$, we then have 
$0=F_\mu'(\phi)=\MA(\phi)-\mu$ as required.

Now, there are at least three problems with this approach:
\begin{itemize}
\item[(1)]
  There is a priori no reason why a maximizer should exist in $\PSH^0(\Lan)$.
  We resolve this by introducing a larger space $\PSH(\Lan)$ with 
  suitable compactness properties and find a maximizer there.
\item[(2)]
  Granted the existence of a maximizer $\phi\in\PSH(\Lan)$, 
  we are maximizing over a convex set rather than 
  a vector space, so there is no reason why $F_\mu'(\phi)=0$.
  Compare maximizing the function $f(x)=x^2$ on the real
  interval $[-1,1]$: 
  the maximum is not at a critical point.
\item[(3)]
  In the end we want to show that---after all---the maximizer is continuous,
  that is, $\phi\in C^0(\Lan)$.
\end{itemize}
We shall discuss how to address~(1) and~(2) in the next two sections.
The continuity result in~(3) requires a priori capacity estimates due to 
Ko{\l}odziej, and will not be discussed in these notes.
%
%
%
%
%
%
\section{Singular semipositive metrics}\label{S106}
Plurisubharmonic (psh) functions are among the \emph{objets souples}
(soft objects) in complex analysis according to P.~Lelong~\cite{Lel85}.
This is reflected in certain useful compactness properties.
The global analogues of psh functions are semipositive singular
metrics on holomorphic line bundles.
Here ``singular'' means that vectors may have infinite length.

\begin{Thm}\label{T104}
  Let $K$ be either $\C$ or a discretely valued field of residue 
  characteristic zero, and let $(X,L)$ be a smooth projective polarized
  variety over $K$.
  Then there exists a unique class $\PSH(\Lan)$,
  the set of singular semipositive metrics, with the following properties:
  \begin{itemize}
  \item
    $\PSH(\Lan)$ is a convex set which is closed under maxima and addition 
    of constants;
  \item
    $\PSH(\Lan)\cap C^0(\Lan)=\PSH^0(\Lan)$;
  \item
    if $s_i$, $1\le i\le p$, are nonzero global sections of $mL$ for some
    $m\ge1$, then $\phi:=\frac1{m}\max_i\log|s_i|\in\PSH(\Lan)$;
    further, $\phi$ is continuous iff the sections $s_i$ have no common zero.
  \item
    if $(\phi_j)$ is an arbitrary family in $\PSH(\Lan)$ that is uniformly 
    bounded from above, then the usc regularization of $\sup_j\phi_j$
    belongs to $\PSH(\Lan)$;
  \item
    if $(\phi_j)$ is a decreasing net in $\PSH(\Lan)$, then either 
    $\phi_j\to-\infty$ uniformly on $\Xan$, or $\phi_j\to\phi$
    pointwise on $\Xan$ for some $\phi\in\PSH(\Lan)$;
  \item
    \textbf{Regularization}: for every $\phi\in\PSH(\Lan)$ there exists a decreasing 
    sequence $(\phi_m)_{m=1}^\infty$ of smooth/model metrics
    such that $\phi_m$  converges pointwise to $\phi$ on $\Xan$ as $m\to\infty$;
  \item
    \textbf{Compactness}: the space $\PSH(\Lan)/\R$ is compact.
  \end{itemize}
\end{Thm}
To make sense of the compactness statement we need to specify the 
topology on $\PSH(\Lan)$. In the complex case, one usually fixes a 
volume form $\mu$ on $\Xan$ and takes the topology
induced by the $L^1$-norm:
$\|\phi-\psi\|=\int_{\Xan}|\phi-\psi|\mu$.
In the non-Archimedean case, there is typically no volume form on $\Xan$.
Instead, we say that a net $(\phi_j)_j$ in $\PSH(\Lan)$ converges to 
$\phi$ if $\lim_j\sup_{\Delta_\cX}|\phi_j-\phi|=0$ for every SNC model $\cX$.
Implicit in this definition is that the restriction to $\Delta_\cX$ of every
singular metric in $\PSH(\Lan)$ is continuous: see
Theorem~\ref{T102} below.

In the complex case, one typically defines $\PSH(\Lan)$ as the 
set of usc singular metrics $\phi$ that are locally represented by 
$L^1$ functions and whose curvature current $dd^c\phi$
(computed in the sense of distributions) is a positive closed current.
Thus $\phi$ is locally given as the sum of a smooth function and a
psh function.
Most of the statements above then follow from basic facts about
plurisubharmonic functions in $\C^n$. 
The regularization result is the most difficult. On $\C^n$ it is easy to
regularize using convolutions. With some care, one can in the global (projective) case
glue together local regularizations to obtain a global one. 
See~\cite{Dem92} for a general result and~\cite{BK} 
for a relatively simple argument applicable in our setting.

In the non-Archimedean case, we are not aware of any workable a priori 
definition of $\PSH(\Lan)$. Chambert-Loir and Ducros~\cite{CD12} 
have a notion of
forms and currents on Berkovich spaces, but it is unclear if it gives the right
objects for the purposes of the theorem above. 
Instead, we prove the following result:
\begin{Thm}\label{T102}
  For any SNC model $\cX$, the restriction of the dual 
  complex $\Delta_\cX\subset\Xan$ of the set of model metrics
  on $\Lan$ forms an equicontinuous family.
\end{Thm}
This is proved using a rather subtle argument, involving intersection numbers
on toroidal models dominating $\cX$. It would be interesting to have a different proof.
At any rate, Theorem~\ref{T102} allows us to define $\PSH(\Lan)$ as the set of 
usc singular metrics $\phi$ satisfying, for every sufficiently large SNC model $\cX$,
\begin{itemize}
\item[(i)]
  $(\phi-\phi_0)\circ r_\cX\ge\phi-\phi_0$;
\item[(ii)]
  the restriction of $\phi$ to $\Delta_\cX$ is a 
  uniform limits of a sequence $\phi_m|_{\Delta_\cX}$, where
  each $\phi_m$ is a semipositive model metric.
\end{itemize}
Here $\phi_0$ is a fixed model metric, determined by some model dominated by $\cX$. 
The map $r_\cX:\Xan\to\Delta_\cX\subset\Xan$ is a natural retraction.
Since $\phi$ is usc, condition~(i) implies that $\phi=\phi_0+\lim_\cX(\phi-\phi_0)\circ r_\cX$,
so that $\phi$ is determined by its restrictions to all dual complexes.

With this definition, the compactness of $\PSH(\Lan)/\R$ follows from Theorem~\ref{T102}
and Ascoli's theorem. Regularization, however, is quite difficult to show. We are not aware of
any procedure that would replace convolution in the complex case. Instead we use algebraic 
geometry. Here is an outline of the proof.

Fix $\phi\in\PSH(\Lan)$. For any SNC model $\cX$, $\phi$ naturally 
induces a model metric $\phi_\cX$. The semipositivity of $\phi$ implies
that the net $(\phi_\cX)_\cX$, indexed by the collection of 
(isomorphism classes of) SNC models decreases to $\phi$.
Unfortunately, except in the curve case $n=1$, 
$\phi_\cX$ has no reason to be semipositive; this reflects the fact
that the pushforward of a nef line bundle may fail to be nef.
We address this by defining $\psi_\cX$ as 
the supremum of all semipositive (singular) metrics dominated by $\phi_\cX$.
We then show that $\psi_\cX$ is continuous and can be \emph{uniformly}
approximated by a sequence $(\psi_{\cX,m})_{m}^\infty$ of semipositive
model metrics. From this data it is not hard to produce a decreasing 
net of semipositive model metrics converging to $\phi$.

Let us say a few words on the construction of the semipositive model
metrics $\phi_{\cX,m}$ since this is a key step in the paper~\cite{siminag}.
For simplicity assume that $L$ is base point free and that 
$\phi_\cX$ is associated to a line bundle $\cL$
(rather than an $\R$-line bundle) on $\cX$. Let $\fa_m$ be the base ideal
of $m\cL$, cut out by the global sections;
it is cosupported on the special fiber $\cX_0$.
The sequence $(\fa_m)_m$ is a \emph{graded sequence} in the sense that $\fa_l\cdot\fa_m\subset\fa_{l+m}$,
Each $\fa_m$ naturally defines a semipositive model metric 
$\psi_{\cX,m}$ on $\Lan$. 
The fact that 
$\psi_{\cX,m}$ converges \emph{uniformly} to $\psi_\cX$ translates into a statement
that the graded sequence $(\fa_m)_m$ is ``almost'' finitely generated. 
This in turn is proved using \emph{multiplier ideals} and ultimately reduces
to the Kodaira vanishing theorem; to apply the latter it is crucial to work
in residue characteristic zero.

The argument above proves that any $\phi\in\PSH(\Lan)$ is the limit
of a decreasing \emph{net} of semipositive model metrics.
When $\phi$ is continuous, the convergence is uniform by 
Dini's Theorem, and we can use the sup-norm to 
extract a decreasing \emph{sequence} of model metrics converging to
$\phi$. In the general case, the \emph{Monge-Amp\`ere capacity} 
developed in~\cite[\S4]{nama} (and modeled on~\cite{BT1,GZ1})
can similarly be used to extract a cenvergent sequence from a net.

%
%
%
%
%
%
\section{Energy}\label{S107}
In the complex case, the (Aubin-Mabuchi) energy functional is defined as follows.
Fix a smooth semipositive reference metric $\phi_0$ and set
\begin{equation}\label{eq:energy}
  E(\phi)
  :=\frac1{n+1}\sum_{j=0}^n
  \int_{\Xan}(\phi-\phi_0)(dd^c\phi)^j\wedge(dd^c\phi_0)^{n-j}.
\end{equation}
for any smooth metric $\phi$. 
Here $(dd^c\phi)^j\wedge(dd^c\phi_0)^{n-j}$ is a 
\emph{mixed Monge-Amp\`ere measure}. It is a positive measure
if $\phi$ is semipositive.

In the non-Archimedean case, mixed Monge-Amp\`ere 
measures can be defined using intersection
theory when $\phi$ and $\phi_0$ are model metrics, and the
energy of $\phi$ is then defined exactly as above.

For two smooth/model metrics $\phi$, $\psi$ we have
\begin{equation}\label{eq:cocycle}
  E(\phi)-E(\psi)
  =\frac1{n+1}\sum_{j=0}^n\int_{\Xan}(\phi-\psi)(dd^c\phi)^j\wedge(dd^c\psi)^{n-j}.
\end{equation}
This is proved using integration by parts in the complex case and
follows from basic intersection theory in the non-Archimedean case.

We can draw two main conclusions from~\eqref{eq:cocycle}.
First, the derivative of the energy functional is the Monge-Amp\`ere
operator, in the sense that
\begin{equation}\label{e101}
  \frac{d}{dt}E(\phi+tf)\bigg|_{t=0}=\int_{\Xan}f\MA(\phi)
\end{equation}
for a smooth/model metric $\phi$ on $\Lan$ and a smooth/model 
function $f$ on $\Xan$.

Second, $E(\psi)\ge E(\phi)$ when $\psi\ge\phi$ are semipositive.
It then makes sense to set
\begin{equation*}
  E(\phi):=\inf \{E(\psi)\mid \psi\ge\phi,\ \text{$\psi$ a
    semipositive smooth/model metric on $\Lan$}\}.
\end{equation*}
for any singular semipositive metric $\phi\in\PSH(\Lan)$.
The resulting functional 
\begin{equation*}
  E:\PSH(\Lan)\to[-\infty,\infty)
\end{equation*}
has many good properties: $E$ is concave, monotonous, and satisfies
$E(\phi+c)=E(\phi)+c$ for $c\in\R$. 
Further, $E$ is usc and continuous along decreasing nets.

The energy functional singles out a class $\cE^1(\Lan)$ of metrics
with \emph{finite energy}, $E(\phi)>-\infty$. This class has good
properties. In particular, one can (with some effort) 
define mixed Monge-Amp\`ere
measures $(dd^c\phi)^j\wedge(dd^c\psi)^{n-j}$ for
$\phi,\psi\in\cE^1(\Lan)$, and~\eqref{eq:energy} continues to hold.

Let us now go back to the variational approach to solving the
Monge-Amp\`ere equation. Fix a  positive measure $\mu$ 
on $\Xan$ of mass $(L^n)$. In the complex case we assume
$\mu$ is absolutely continuous with respect to Lebesgue measure, with 
density in $L^p$ for some $p>1$. In the non-Archimedean case
we assume that $\mu$ is supported on some dual complex. 
In both cases, one can show that the functional
$\phi\to\int(\phi-\phi_0)\,\mu$ is (finite and) continuous
on $\PSH(\Lan)$, where $\phi_0$ is the same reference metric as
in~\eqref{eq:energy}. Thus the functional 
$F_\mu\colon\PSH(\Lan)\to[-\infty,\infty)$ defined by 
\begin{equation*}
  F_\mu(\phi):=E(\phi)-\int(\phi-\phi_0)\mu
\end{equation*}
is upper semicontinuous. It follows from~\eqref{eq:cocycle} that 
$F_\mu$ does not depend on the choice of reference metric $\phi_0$.
We also have $F_\mu(\phi+c)=F_\mu(\phi)$ for $\phi\in\PSH(\Lan)$, 
$c\in\R$. Thus $F_\mu$ descends to an usc functional on the 
quotient space $\PSH(\Lan)/\R$. 
By Theorem~\ref{T104}, the latter space is compact, so we can find
$\phi\in\PSH(\Lan)$ maximizing $F_\mu$.
It is clear that $\phi\in\cE^1(\Lan)$, so the mixed Monge-Amp\`ere
measures of $\phi$ and $\phi_0$ are well defined.
However, equation~\eqref{e101} no longer makes sense, since
there is no reason for the metric $\phi+tf$ to be semipositive for
$t\ne0$. Therefore, it is not clear that $\MA(\phi)=\mu$,
as desired.
In the next section, we explain how to get around this problem.
%
%
%
%
%
%
\section{Envelopes, differentiability and orthogonality}\label{S108}
We define the \emph{psh envelope} of 
a (possibly singular) metric $\psi$ on $\Lan$ by 
\begin{equation*}
  P(\psi):=\sup\{\phi\in\PSH(\Lan)\mid \phi\le\psi\}^*.
\end{equation*}
As before, $\phi^*$ denotes the usc regularization of a singular metric $\phi$.
In all cases we need to consider, $\psi$ will be the sum of a 
metric in $\cE^1(\Lan)$ and a continuous function on $\Xan$. 
In particular, $\psi$ is usc, $P(\psi)\in\cE^1(\Lan)$ and $P(\psi)\le \psi$.

\begin{figure}
  \includegraphics[width=8cm]{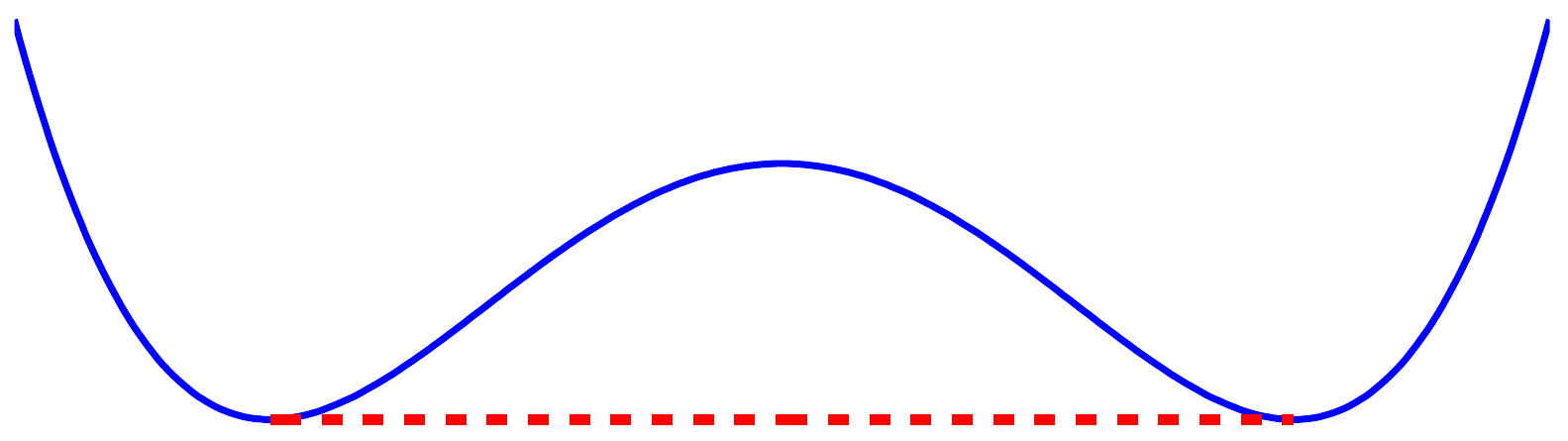}
  \caption{The convex hull $P(f)$ of a continuous function $f$ of one
    variable. Note that $P(f)$ is affine, \ie $P(f)''=0$ where $P(f)\ne f$.}\label{F101}
\end{figure}

This envelope construction was in fact already mentioned at the 
end of~\S\ref{S106} as it plays a key role in the regularization
theorem. The psh envelope is an analogue of the convex hull;
see Figure~\ref{F101}.

The key fact about the psh envelope is that the composition
$E\circ P$ is differentiable and that $(E\circ P)'=E'\circ P$.
More precisely, we have:
\begin{Thm}
  For any $\phi\in\cE^1(\Lan)$ and $f\in C^0(\Xan)$, the function 
  $t\mapsto E(P(\phi+tf))$ is differentiable at $t=0$, with derivative
  $\frac{d}{dt}E(\phi+tf)|_{t=0}=\int f\MA(\phi)$.
\end{Thm}
Granted this result, let us show how to solve the Monge-Amp\`ere equation.
Pick $\phi\in\cE^1(\Lan)$ that maximizes $F_\mu(\phi)=E(\phi)-\int(\phi-\phi_0)\mu$
and consider any  $f\in C^0(\Xan)$. For any $t\in\R$ we have 
\begin{align*}
  E(P(\phi+tf))-\int(\phi+tf-\phi_0)\mu
  &\le  E(P(\phi+tf))-\int(P(\phi+tf)-\phi_0)\mu\\
  &\le E(\phi)-\int(\phi-\phi_0)\mu.
\end{align*}
Since the left hand side is differentiable at $t=0$, the derivative 
must be zero, which amounts to $\int f\MA(\phi)-\int f\mu=0$.
Since $f\in C^0(\Xan)$ was arbitrary, this means that $\MA(\phi)=\mu$,
as desired.

The proof of this differentiability results proceeds by first reducing to the case 
when $\phi$ and $f$ are continuous. A key ingredient is then
\begin{Thm}\label{T101}
  For any continuous metric $\phi$ on $\Lan$ we have
  \begin{equation}\label{eq:orth}
    \int_{\Xan}(\phi-P(\phi))\MA(P(\phi))=0.
  \end{equation}
\end{Thm}
In other words, the Monge-Amp\`ere measure $\MA(P(\phi))$
is supported on the locus $P(\phi)=\phi$. A version of this for
functions of one variable is illustrated in Figure~\ref{F101}.

To prove this result, we can reduce to the case when $\phi$ is a
smooth/model metric. In the complex case, Theorem~\ref{T101}
was proved by Berman and the first author in~\cite{BB} using the 
pluripotential theoretic technique known as ``balayage''. 
In the non-Archimedean setting, Theorem~\ref{T101} is deduced
in~\cite{nama} from the asymptotic orthogonality of Zariski 
decompositions in~\cite{BDPP} and is for this reason called
the \emph{orthogonality property}. 
The assumption in Theorem~\ref{T103} that the variety $X$ be 
defined over a smooth $k$-curve is used exactly in order 
to apply the result from~\cite{BDPP}.

\medskip
The solution to the non-Archimedean Monge-Amp\`ere equation 
$\MA(\phi)=\mu$ can be made slightly more explicit in the case when
the support of $\mu$ is a singleton, $\mu=d_L\,\delta_x$,
where $d_L:=(L^n)$ and $x\in\Xan$ belongs to some dual complex;
such points $x$ are known as \emph{quasimonomial}
or \emph{Abhyankar} points.

The fiber $L^{\mathrm{an}}_x$ of $\Lan$ above $x\in\Xan$ is isomorphic to 
the Berkovich affine line over the complete residue field $\cH(x)$.
Fix any nonzero $y\in L^{\mathrm{an}}_x$ and set 
\begin{equation*}
  \phi_x:=\sup\{\phi\in\PSH(\Lan)\mid\|y\|_\phi\ge 1\}.
\end{equation*}
By~\cite[Prop.8.6]{nama}, $\MA(\phi_x)$ is supported
on $x$, so $\MA(\phi_x)=d_L\,\delta_x$.
It would be interesting to find an example of a divisorial point 
$x\in\Xan$ such that $\phi_x$ is not a model function.
%
%
%
%
%
%
\section{Curves}\label{S109}
The disadvantage of the variational approach to the Monge-Amp\`ere equation
is that it gives very little
control on the solution beyond continuity. Here we shall make 
the solution more concrete in the case of curves; the next section
deals with toric varieties.

Thus assume $X$ is a smooth projective curve over $K$.
In this case, the Monge-Amp\`ere operator (which one would
normally refer to as the Laplacian) is \emph{linear}:
if $\phi_i$ is a metric on $L_i$, $i=1,2$, then 
$\MA(\phi_1+\phi_2)=\MA(\phi_1)+\MA(\phi_2)$.
Furthermore, as we shall see, the Monge-Amp\`ere operator is
naturally defined on \emph{any} singular semipositive metric on $\Lan$,
for an ample line bundle $L$, and we can solve $\MA(\phi)=\mu$ 
for \emph{any} positive measure $\mu$ of mass $\deg L$.

\medskip
Let us first explain this in the complex case; $\Xan$ is then 
a compact Riemann surface. Fix a smooth metric $\phi_0$ on $\Lan$.
The curvature form $\omega_0:=dd^c\phi_0$ is a volume form of
mass $\deg L$. A singular metric $\phi$ on $\Lan$ is then semipositive iff
$\f:=\phi-\phi_0$ is an \emph{$\omega_0$-psh function}, 
that is, a locally integrable function $\f$ that is locally the sum of a
smooth function and a psh function, and such that $\omega_0+dd^c\f$
is a positive measure. 
We then set $\MA(\phi):=\omega_0+dd^c\phi$; this definition does not
depend on the choice of $\phi_0$.

Now we explain how to solve the equation 
$\MA(\phi)=\mu$ for any positive measure $\mu$ of mass 
$d_L:=\deg L$.
Writing $\phi=\phi_0+\f$ as above, we must solve 
$dd^c\f=\mu-\omega_0$, where $\omega_0:=dd^c\phi_0$.
It sufficies to do this when $\mu=d_L\delta_x$ for some 
$x\in\Xan$: indeed,  if we normalize the solution $\f_x$ to
$dd^c\f_x=\mu-d_L\delta_x$ by $\int_{\Xan}\f_x\,\omega_0=0$,
then the function $\f_\mu$ defined by 
$\f_\mu(y):=d_L^{-1}\int_{\Xan}\f_x(y)\,d\mu(x)$ satisfies
$dd^c\f_\mu=\omega_0-\mu$ and is normalized by 
$\int_{\Xan}\f_\mu\,\omega_0=0$.

The function $\f_x$ can be ``physically'' interpreted as the
voltage (suitably normalized) when putting a charge of $+d_L$ 
at the point $x$ and a total
charge of $-d_L$ spread out according to the measure
$\omega_0$. Mathematically, Perron's method describes it as
the supremum of all $\omega_0$-subharmonic functions $\f$ on $\Xan$
satisfying $\int_{\Xan}\f\,\omega_0=0$ and $\f\le d_L\log|z|+O(1)$, 
where $z$ is a local coordinate at $x$.

\medskip
Now we consider the non-Archimedean case. As before, let us
assume that $K$ is a discretely valued field of residue characteristic
zero, even though this is not really necessary in the one-dimensional
case.\footnote{Indeed, Thuillier~\cite{thuillierthesis} systematically develops a potential theory on Berkovich curves in a very general setting.}

The main point is that any Berkovich curve has the structure of a
generalized\footnote{This means that some distances may be infinite.} 
metric graph.
We will not describe this in detail, but here is the idea.  
The dual graph $\Delta_\cX$ of any SNC model $\cX$ is a 
connected, one-dimensional simplicial complex. 
As before, we view it as a subset of $\Xan$.
It carries a natural integral affine structure, inducing a metric.
If $\cX'$ is an SNC model dominating $\cX$ (in the sense that the 
canonical birational map $\cX\dashrightarrow\cX'$ is a morphism)
then $\Delta_\cX$ is a subset of $\Delta_{\cX'}$ and 
the inclusion $\Delta_\cX\to\Delta_{\cX'}$ is an isometry.
There is a also a (deformation) retraction
$r_\cX:\Xan\to\Delta_\cX$, and 
$\Xan\simeq\varprojlim_\cX\Delta_\cX$.
In this way, the metrics on the dual complexes induce a generalized metric on $\Xan$.

The structure of each $\Delta_\cX$ and of $\Xan$ as metric graphs allows
us to define a Laplacian on these spaces, by combining the real Laplacian on 
segments and the combinatorial Laplacian at branch points (and
endpoints).
This Laplacian allows us to understand both semipositive singular
metrics and the Monge-Amp\`ere operator.

Namely, fix a model metric $\phi_0$ on $\Lan$. 
It is represented by a $\Q$-line bundle on some SNC model $\cX^{(0)}$.
The measure $\omega_0:=dd^c\phi_0$ is supported on the vertices of 
$\Delta_{\cX^{(0)}}$.
Now, a singular metric $\phi$ is semipositive iff for every SNC model
$\cX$ dominating $\cX^{(0)}$, the restriction of the function
$\phi-\phi_0$ to $\Delta_\cX$ is a $\omega_0$-subharmonic function
in the sense that
$\Delta((\phi-\phi_0)|_{\Delta_\cX})=\mu_\cX-\omega_0$,
where $\mu_\cX$ is a positive measure on $\Delta_\cX$ of mass $d_L$.
In this case, there further exists a unique measure $\mu$ on $\Xan$
of mass $d_L$ such that $(r_\cX)_*\mu=\mu_\cX$ for all $\cX$,
and we have $\MA(\phi)=\mu$.

To solve the equation $\MA(\phi)=\mu$ for a positive measure $\mu$ of
mass $d_L$, it suffices by linearity to treat the case $\mu=d_L\delta_x$
for a point $x\in\Xan$. In this case, the function $\f:=\phi-\phi_0$ will
be locally constant outside the convex hull of
$\{x\}\cup\Delta_{\cX^{(0)}}$.
The latter is essentially a finite metric graph on which we need to 
find a function whose Laplacian is equal to $d_L\delta_x-\omega_0$.
This can be done in a quite elementary way.

\medskip
An interesting example of semipositive metrics, both in the complex
and non-Archimedean case, comes from dynamics~\cite{Zha95a}.
Suppose $f:(X,L)\self$ is an polarized endomorphism of degree
$\lambda>1$. In other words, $f:X\to X$ is an endomorphism and 
$f^*L$ is linearly equivalent to $\lambda L$.
Then there exists a unique \emph{canonical metric} $\phi_{\can}$ on $\Lan$,
satisfying $f^*\phi=\lambda\phi$. This metric is continuous and
semipositive but usually not a model metric.

As a special case, suppose $X$ is an elliptic curve and that $f$ is
the map given by multiplication by $\lambda$. 
In the complex case, $\Xan\simeq\C/\Lambda$ is a torus 
and $\mu_{\can}:=\MA(\phi_{\can})$ 
is given by a multiple of Haar measure on $\Xan$.
In the non-Archimedean case, there are two possibilities. 
If $X$ has good reduction over $\Spec R$, then $\mu_{\can}$ is 
a point mass. Otherwise, $\mu_{\can}$ is proportional to
Lebesgue measure on the \emph{skeleton} $\mathrm{Sk}(\Xan)$,
a subset homeomorphic to a circle. 
A similar description of the measure $\mu_{\can}$ in the case
of higher-dimensional abelian varieties is given in~\cite{Gub10}.
%
%
%
%
%
%
\section{Toric varieties}\label{S110}
For general facts about toric varieties, 
see~\cite{FultonToric,KKMS,BPS}. In this section we briefly describe 
how the complex and non-Archimedean points of view
elegantly come together in the toric setting and translate into
statements about convex functions and the real Monge-Amp\`ere operator.
As before, we only consider the non-Archimedean field
$K=k\lau{t}$ with $\charac k=0$;
however, most of what 
we say here should be true in a more general context: 
see~\cite{Gub13a}.

Let $M\simeq\Z^n$ be a free abelian group, $N$ its dual, 
and let $T=\Spec K[M]$ be the corresponding split $K$-torus. 
A polarized toric variety $(X,L)$ is then determined by a 
rational polytope $\Delta\subset M_\R$.
The variety $X$ is described by the normal fan to $\Delta$ in $N_\R$
and the points of $M\cap\Delta$ are in 1-1 correspondence 
with equivariant sections of $L$; we write $\chi^u$ for the 
section of $L$ associated to $u\in M$.
This description is completely general and holds over any field 
as well as over $\Z$. 

There is also a ``tropical'' space $\Xtrop$ associated to $X$.
As a topological space, it is compact and contains $N_\R$ as an 
open dense subset.\footnote{In our setting, 
$\Xtrop$ can be identified with the 
(moment) polytope $\Delta$ in such a way that $N_\R$
corresponds to the interior of $\Delta$, but this identification 
does not preserve the affine structure on $N_\R$.}
For any valued field $K$, there is a tropicalization map 
$\trop:\Xan\to\Xtrop$, where $\Xan$ refers to the analytification with 
respect to the norm on $K$.
The inverse image of $N_\R$ is the torus $\Tan$.

There is a natural correspondence between equivariant metrics on 
$\Lan$ and functions on $N_\R$. Let
$\phi$ is an equivariant metric on $\Lan$. 
For every $u\in M$, $\chi^u$, is a nonvanishing section of $L$ on $T$
so $\phi-\log|\chi^u|$ defines a function on 
$\Tan$ that is constant on the fibers of the tropicalization map.
In particular, picking $u=0$, we can write 
\begin{equation}\label{e102}
  \phi-\log|\chi^0|=g\circ\trop
\end{equation}
for some function $g$ on $N_\R$.
Conversely, given a function $g$ on $N_\R$,~\eqref{e102} defines
an equivariant metric on the restriction of $\Lan$ to $\Tan$.

We now go from the torus $T$ to the polarized variety $(X,L)$.
After replacing $L$ by a multiple, we may assume that all the vertices
of $\Delta$ belong to $M$. Set
\begin{equation*}
  \phi_\Delta:=\max_{u\in\Delta}\log|\chi^u|.
\end{equation*}
This is a semipositive, equivariant model metric on $\Lan$.
Its restriction to $\Tan$ corresponds to the 
homogeneous, nonnegative, convex function
\begin{equation*}
  g_\Delta:=\max_{u\in\Delta}u
\end{equation*}
on $N_\R$. 
In general, an equivariant singular metric $\phi$ on $\Lan$ corresponds 
to a convex function $g$ on $N_\R$ such that $g\le g_\Delta+O(1)$.
It is bounded iff $g-g_\Delta$ is bounded on $N_\R$.

\medskip
The real Monge-Amp\`ere measure of any convex function $g$ on $N_\R$ is a well-defined positive measure $\MA_\R(g)$ on $N_\R$ (see e.g. \cite{RT}).
When $g=g_\Delta+O(1)$, its total mass is given by 
\begin{equation*}
  \int_{N_\R}\MA_\R(g)=\Vol(\Delta)=\frac{(L^n)}{n!},
\end{equation*}
where the last equality follows from~\cite[p.111]{FultonToric}.

\medskip
We now wish to relate the real Monge-Amp\`ere measure of $g$ 
and the Monge-Amp\`ere measure of the corresponding 
semipositive metric $\phi$ on $\Lan$.

First consider the non-Archimedean case, in which 
there is a natural embedding $j:N_\R\to\Tan\subset\Xan$ 
given by monomial valuations that sends
$v\in N_\R$ to the norm
\begin{equation*}
  \sum_{u\in M} a_uu \in K[M]
  \mapsto\max_{u\in M}\{|a_u|\exp(-\langle u,v\rangle)\}.
\end{equation*}
In particular, $j(0)=x_G$, the Gauss point of the open $T$-orbit. 

If $g$ is a convex function on $N_\R$ with $g=g_\Delta+O(1)$, 
and if $\phi$ is the corresponding continuous semipositive metric on $L$, then~\cite[Theorem~4.7.4]{BPS} asserts that
\begin{equation*}
  \MA(\phi)=n!\,j_*\MA_\R(g). 
\end{equation*}
For a compactly supported positive measure $\nu$ on $N_\R$
of mass $(L^n)$, 
solving the Monge-Amp\`ere equation $\MA(\phi)=j_*(\nu)$ 
therefore amounts to solving the real 
Monge-Amp\`ere equation $\MA_\R(g)=\nu/n!$.
This can be done explicitly when $\nu$ is a point mass,
say supported at $v_0\in N_\R$. Indeed, the function 
$g_{v_0}:N\to\R$ defined by $g=g_\Delta(\cdot-v_0)$ is convex
and satisfies $g=g_\Delta+O(1)$. Further, for every point 
$v\ne v_0$ there exists a line segment in $N_\R$ containing $v$ 
in its interior and on which $g$ is affine. This implies
that $\MA_\R(g)$ is supported at $v_0$. As a a consequence,
the corresponding continuous metric $\phi$ on $\Lan$
satisfies $\MA_\R(\phi)=(L^n)\delta_{j(u_0)}$.

This solution can be shown to tie in well with the construction at the 
end of~\S\ref{S108}, but is of course much more explicit. For example,
when $u_0\in N_\Q$, so that $j(u_0)\in\Xan$ is divisorial,
the function $g_{u_0}$ is $\Q$-piecewise linear so that the corresponding
metric $\phi$ is a model metric.

\medskip
Finally we consider the complex case. In this case we cannot embed $N_\R$
in $\Tan$. However, the preimage of any point 
$v\in N_\R$ under the tropicalization is a real torus of dimension 
$n$ in $\Tan$ on which the multiplicative group $(S^1)^n$ acts
transitively. 
To any compactly supported positive 
measure $\nu$ on $N_\R$ of mass $(L^n)/n!$ we can therefore
associate a unique measure $\mu$ on $\Tan$, still denoted 
$\mu:=j_*\nu$, that is invariant under
the action of $(S^1)^n$ and satisfies $\trop_*\mu=\nu$. 

If $\phi$ is an equivariant semipositive metric on $\Lan$,
corresponding to a convex function $g$ on $N_\R$, we then have
\begin{equation*}
  \MA(\phi)=n!\,j_*\MA_\R(g). 
\end{equation*}
For $(S^1)^n$-invariant measures $\mu$ on $\Lan$ of mass $(L^n)$,
solving the complex Monge-Amp\`ere equation $\MA(\phi)=\mu$
thus reduces to solving the real Monge-Amp\`ere equation 
$\MA_\R(g)=\frac1{n!}\trop_*\mu$. 
%
%
%
%
%
%
\section{Outlook}\label{S111}
In this final section we indicate some possible extensions of our work
and make a few general remarks.

\medskip
First of all, it would be nice to have a \emph{local} theory for
semipositive singular metrics. Indeed, while the global approach
in~\cite{siminag,nama} serves works well for the Calabi-Yau problem,
it has some unsatisfactory features. For example, it is not completely trivial to
prove that the Monge-Amp\`ere operator is local in the sense 
that if $\phi_1,\phi_2$ are two (say) continuous semipositive metrics 
that agree on an open subset $U\subset\Xan$, 
then $\MA(\phi_1)=\MA(\phi_2)$ on $U$. We prove this in~\cite{nama}
using the Monge-Amp\`ere capacity. 
Still, it would be desirable to say that the restriction of a semipositive
metric to (say) an open subset of $\Xan$ remains semipositive!

In contrast, in the complex case, the classical approach is local in nature. 
Namely, one first defines and studies psh functions 
on open subsets of $\C^n$ and then defines
singular semipositive metrics as global analogues. 
By construction, the Monge-Amp\`ere operator is a 
local (differential) 
operator.\footnote{However, one also needs to verify that the Monge-Amp\`ere operator is local for the \emph{plurifine}
topology. This is nontrivial in both the complex and non-Archimedean case.}

In a general non-Archimedean setting, 
Chambert-Loir and Ducros~\cite{CD12} (see also~\cite{Gub13b,GK14})
define psh functions as continuous functions $\f$ such that 
$d'd''\f$ is a positive closed current (in their sense), for suitable 
operators $d'$, $d''$ analogous to their complex counterparts and
modeled on notions due to Lagerberg~\cite{Lag12}.
While this leads to a very nice theory, that moreover works 
for general Berkovich spaces, the crucial compactness
and regularization results are so far missing.
At any rate, the tropical charts used in~\cite{CD12} may be a 
good substitute for dual complexes of SNC models.

\medskip
Going back to the projective setting, there are several open questions
and possible extensions, even in the case of a discretely valued
ground field of residue characteristic zero.

First, when solving the Monge-Amp\`ere 
equation, we needed to assume that the variety $X$ was obtained by base change from a variety over a $k$-curve.  This assumption was made in 
order to use the orthogonality result in~\cite{BDPP}, but is presumably redundant.

Second, one should be able to solve the Monge-Amp\`ere equation 
$\MA(\phi)=\mu$ for more general measures $\mu$. In the complex setting,
this is done in~\cite{GZ2,Dinew} for non-pluripolar measures $\mu$.
The analogous result should be valid in the non-Archimedean setting, too, although some countability issues seem to require careful attention.
Having such a general result would allow for a nice Legendre duality, as
explored in~\cite{BBGZ,Berm13} in the complex case. 

Third, one could try to get more specific information about the 
solution. We already mentioned at the end of~\S\ref{S108} that we
don't know whether the solution to the equation $\MA(\phi)=d_L\,\mu_x$
is a model function for $x$ a divisorial point (and $d_L=(L^n)$).
In a different direction, one could consider the case when $X$ is a 
Calabi-Yau variety, in the sense that $K_X\simeq\cO_X$. 
Then there exists a canonical subset $\mathrm{Sk(X)}\subset\Xan$,
the \emph{Kontsevich-Soibelman skeleton}, 
see~\cite{KoSo,mustata-nicaise,nicaise-xu}.
It is a subcomplex of the dual complex of any SNC model 
and comes equipped
with an integral affine structure, inducing a volume form on each 
face. One can solve $\MA(\phi)=\mu$, for linear combinations
of these volume forms, viewed as measures on $\Xan$.
Can we say anything concrete about the solution $\phi$,
as in the case of maximally degenerate 
abelian varieties considered in~\cite{Liu11}?

\medskip
It would obviously be interesting to work
over other types of non-Archimedean fields, such as $\Q_p$.
Here there are several challenges. First, we systematically use
SNC models, which are only known to exist in residue 
characteristic zero (except in low dimensions). 
It is possible that the tool of SNC models can, with some additional effort, 
be replaced by alterations, tropical charts or other methods. 
However, we also crucially
use the assumption of residue characteristic zero when applying
the vanishing theorems that underly the regularization theorem for 
singular semipositive metrics. Here some new ideas are needed.

A simpler situation to handle is that of a \emph{trivially} valued
field. This is explored in~\cite{trivval} and can be briefly explained
as follows. Let $k$ be any field of characteristic zero, equipped with
the trivial norm. Let $(X,L)$ be a polarized variety over $k$.
In this setting, the notion of model metrics and model functions seemingly
does not take us very far, as the only model of $X$ is $X$ itself! 
Instead, the idea is to use a non-Archimedean 
field extension. Set $K=k\lau{t}$, $X_K:=X\otimes_kK$ etc.
The multiplicative group $G:=\G_{m,k}$ acts on $X_K^{\mathrm{an}}$
and $\Xan$ can be identified with the set of $G$-equivariant 
points in $X_K^{\mathrm{an}}$. Similarly, singular semipositive
metrics on $\Lan$ are defined as $G$-invariant singular 
semipositive metrics on $L_K^{\mathrm{an}}$.
In this way, the main results about $\PSH(\Lan)$ follow from
the corresponding results about $\PSH(L_K^{\mathrm{an}})$
and the same is true for the solution of the Monge-Amp\`ere equation.

A primary motivation for studying the trivially valued case, at least in
the case $k=\C$, is that 
the space of singular semipositive metrics on $\Lan$
naturally sits ``at the boundary'' of the space of 
positive (K\"ahler) metrics on the holomorphic line bundle $L$.
As such, it can be used to study questions on $K$-stability and may
be useful for the study of the existence of constant scalar curvature metrics,
see~\cite{kunif1,kunif2}.
A different scenario where a complex situation degenerates to
a non-Archimedean one occurs in~\cite{amoebae}.

\medskip
In yet another direction, one could try to consider line bundles that are not necessarily ample, but rather big and nef, or simply big. 
In the complex case this was done in~\cite{EGZ,BEGZ}. 
One motivation for such a generalization is that it is invariant under birational maps and would hence allow us to study singular varieties.

Finally, it would be interesting to have transcendental analogues. Indeed,
it the complex case, one often starts with a K\"ahler manifold $X$ together 
with a K\"ahler class $\omega$, rather than a polarized pair $(X,\omega)$.
A notion of K\"ahler metric is proposed in~\cite{KoTs,YuGromov}, but
it is not clear whether or not this plays the role of a (possibly)
transcendental K\"ahler metric.
%
%
%
%
%
%

\end{document}